\documentclass[twoside]{article}
\usepackage{amsmath,amssymb,amsthm,latexsym,bm,indentfirst,wasysym}
\usepackage{graphicx,epsfig,amscd,mathrsfs,multirow,bm}
\usepackage{cite}
\usepackage{caption}
\usepackage{color}

\input xy
\xyoption{all}
\input scrload.tex

\allowdisplaybreaks[4]

\catcode`@=11
\long\def\@makefntext#1{\noindent #1}
\newskip\tabcentering \tabcentering=1000pt plus 1000pt minus 1000pt
\def\MCH#1#2{\setbox0=\hbox{\raise#1\hbox{#2}}\smash{\box0}}

\def\@evenfoot{}\def\@oddfoot{}

\def\@evenhead{\hbox to\textwidth{\small\rm\thepage \hfill
{\it Qi LIU, Linlin FU and Yongjin LI}}} 

\def\@oddhead{\hbox to \textwidth{\small{\it
On the stability of orthogonal additivity in $\beta$-homogeneous  $F$-spaces 
} \hfill\thepage}}   

\def\SUB#1{\vskip .2in\leftline{\large\bf #1}\vskip .1in}

\catcode`@=12

\def\bc{\begin{center}}
\def\ec{\end{center}}
\def\no{\noindent}
\def\hang{\hangindent\parindent}
\def\textindent#1{\indent\llap{\qquad #1\ \ \enspace}\ignorespaces}
\def\ref{\par\hang\textindent}

\begin{document}


\abovedisplayskip=6pt plus 1pt minus 1pt \belowdisplayskip=6pt
plus 1pt minus 1pt
\thispagestyle{empty} \vspace*{-1.0truecm} \noindent
\vskip 10mm \bc{\Large\bf On The Stability Of Orthogonal Additivity In $\beta$-Homogeneous  $F$-Spaces    
\footnotetext{\footnotesize
Supported by the National Natural Science Foundation of P. R. China (Nos. 11971493 and 12071491). Linlin Fu was supported by the Fundamental Research Funds for the Central Universities, Sun Yat-sen University (Grants 2021qntd21).\\
* Corresponding author\\  
E-mail address: liuq325@mail2.sysu.edu.cn(Qi LIU); fullin3@mail.sysu.edu.cn(Linlin FU);\\ 
stslyj@mail.sysu.edu.cn(Yongjin LI)
} } \ec  

\vskip 5mm
\bc{\bf Qi Liu,\ \ \ Linlin Fu$^{*}$,\ \ \ Yongjin Li}\\  
{\small\it Department of Mathematics,
Sun Yat-sen University, Guangzhou 510275, China
}\ec   

\vskip 1 mm

{\narrower\noindent{\small {\small\bf Abstract}\ \ In this paper, we study the stability of the orthogonal equation, which
is closely related to the results  by  Włodzimierz Fechner and Justyna Sikorska in 2010. There are some differences that we consider the  target space with the  $\beta$-homogeneous norm and quasi-norm. Overcoming the $\beta$-homogeneous norm and quasi-norm bottlenecks, we get some new results.

\vspace{1mm}\baselineskip 12pt

\no{\small\bf Keywords} \ \ Hyers-Ulam stability;  $\beta$-homogeneous  $F$-spaces; quasi-Banach space; orthogonal additivity. 

\par

\vspace{2mm}

\no{\small\bf MR(2020) Subject Classification\ \ {\rm 39B55, 39B82.
}} 

}}

\baselineskip 15pt

\SUB{1. Introduction}
The stability problem of functional equations originated from a problem raised by Ulam \cite{1} on the stability of group homomorphisms in 1940:

Given a group $G$ and a metric group $G^{\prime}$ with metric $\rho(\cdot,\cdot)$. Given $\varepsilon>0$, does there exist a $\delta>0$ such that if $f:G\rightarrow G^{\prime}$ satisfies $\rho(f(xy),f(x)f(y))<\delta$  for all $x,y\in G$, then a homomorphism $h:G\rightarrow G^{\prime}$ exists with $\rho(f(x),h(x))<\varepsilon$ for
all $x\in G$ ?

Hyers \cite{2} gave a first affirmative partial answer to the question of Ulam for additive mappings on Banach spaces. 
Firstly, for the additive mappings, the Hyers theorem generalized form was solved by Aoki \cite{3},  and further, for linear mapping, it was generalized by Rassias \cite{4} taking an unbounded Cauchy difference in consideration. 
A generalization of the Rassias theorem was obtained by Găvruta \cite{5} by replacing the unbounded Cauchy difference by a general control function. Because of their breakthrough achievements, the stability of functional equations has been widely studied by mathematicians.

Although various studies on stability have been successfully conducted, there are not many corresponding stability results due to the non-linear structure of the infinite-dimensional $F$-space.
The nonlinear structure of $F$-space plays an important role in functional analysis and other mathematical fields. 
The $L^p([0,1])$ for $0<p<1$ equipped with the metric $d(f,g)=\int \vert f(x)-g(x) \vert^pdx$ 
is an example of an  $F$-space  but not a Banach space.
Besides these, for $F$-spaces, and we recommend readers to read the literature \cite{6,7}.

{\bf Definition 1.1}\label{dt1}
{\sl Consider $X$ be a linear space. A non-negative valued function $\|\cdot\|$ achieves 
	an $F$-norm if satisfies the following conditions: 
	
	{\rm (1)} $\|x\|= 0$ if and only if $x=0$;
			
	{\rm (2)} $\|\lambda x\|=\|x\|$ for all  $\lambda$,  $|\lambda|=1$;
		
	{\rm (3)} $\|x+y\|\leq\|x\|+\|y\|$ for all $x,y\in X$;
			
	{\rm (4)} $\|\lambda_{n}x\|\rightarrow 0$ provided  $\lambda_{n}\rightarrow 0$;
		
	{\rm (5)} $\|\lambda x_{n}\|\rightarrow 0$ provided  $x_{n}\rightarrow 0$;
		
	{\rm (6)} $\|\lambda_{n} x_{n}\|\rightarrow 0$ provided  $\lambda_{n}\rightarrow 0, x_{n}\rightarrow 0$.

{\it Then $(X,\|\cdot\|)$ is called an $F^\ast$-space. An F-space is a complete $F^\ast$-space}.

An $F$-norm  is called $\beta$-homogeneous $(\beta>0)$ if $\|tx\|=|t|^{\beta}\|x\|$ for all $x\in X$ and all $t\in \mathcal{C}$ (see \cite{8,9}).
}

If a quasi-norm  is $p$-subadditive, then it is called $p$-norm $(0 < p < 1)$.
In other words, if it satisfies
$$\Vert x+y\Vert^p\leq \Vert x\Vert^p+\Vert y\Vert^p,~~x,y\in X.$$
We note that the $p$-subadditive quasi-norm  $\Vert \cdot\Vert $ induces an $F$ norm.
We refer the reader to \cite{6} and  \cite{10} for background on it.

{\bf Definition 1.2}\rm{\cite{11}}
{\sl
 A quasi–norm on $\Vert \cdot\Vert$ on vector space $X$ over a field $K$ $(
 \mathbb{R})$ is a map $X\longrightarrow[0, \infty)$ with the following properties:

 {\rm (1)}
$\|x\|=0$  if and only if $x=0$;

{\rm (2)}
$\|a x\|=|a|\|x\|,~~ a \in \mathbb{R}, x \in X$;

{\rm (3)}
$\|x+y\| \leq C(\|x\|+\|y\|),~~ x, y \in X$.

where $C \geq 1$ is a constant independent of $x, y \in X$. The smallest $C$ for which
$(3)$ holds is called the quasi-norm constant of $(X,\Vert \cdot\Vert)$.}

It is vital to emphasize the well-known theorem in nonlocally convex theory, that is,
Aoki–Rolewicz theorem \cite{8}, which asserts that for some $0 < p\leq 1$, every quasi-norm  
admits  an equivalent $p$-norm.

Various more results for the stability of  functional equations in quasi-Banach spaces
can be seen in  \cite{12,13}. However, the results are more interesting and meaningful 
when orthogonality is taken into account.

The notion of orthogonality goes a long way back in time and various extensions have been introduced over the last decades. In particular, proposing the notion of orthogonality in normed linear spaces has been the object of extensive efforts of many mathematicians.

We recall two orthogonality types introduced in normed linear spaces. In 1945 James \cite{14} introduced the so-called isosceles orthogonality as follows:
$$x\perp_I y ~if~ and ~only ~if~ \Vert x+y\Vert=\Vert x-y\Vert.$$
Taking into account the classical Pythagorean theorem, one can define the orthogonal relation in a normed space $(X, \Vert \cdot \Vert)$:

$$x\perp_P y ~if~ and ~only ~if~  \|x-y\|^{2}=\|x\|^{2}+\|y\|^{2}.$$
Some other known orthogonalities in normed linear spaces can be found in
\cite{15,16,17,18} and references therein.

When orthogonality is the general orthogonality on the inner product space, G. Pinsker \cite{19} characterizes the orthogonality additive functional on the inner product space. K. Sundaresan \cite{20} extends this result to arbitrary Banach Spaces with Birkhoff-James orthogonality.  The following orthogonal Cauchy functional equation was first investigated 
by S. Gudder and D. Strawther \cite{21}
$$
f(x+y)=f(x)+f(y), \quad x \perp y,
$$
where $\perp$ is an abstract orthogonality relation and it
will paly a crucial role in orthogonal stability.

In addition to the different definitions of orthogonality in normed space, we can also give some axiomatic definitions of such relations in linear space.  We show the following standard definition by R\"{a}tz \cite{22}:

{\bf Definition 1.3}
{\sl
	Let $X$ be a real linear space with $\operatorname{dim}X \geq  2$ and let $\perp$ be a binary relation 
	on $X$ such that
	
	{\rm (1)} $x\perp 0$ and $0\perp x$ for all $x\in X$;
	
	{\rm (2)}  if $x,y\in X \backslash \{0\}$  and $x\perp y$, then $x$ and $y$ are linearly independent;
	
	{\rm (3)}	 if $x,y\in X$ and $x\perp y$, then for all $\alpha,\beta \in \mathbb{R}$ we have $\alpha x\perp \beta y$;
	
	{\rm (4)}	  for any two-dimensional subspace $P$ of $X$ and for every $x \in  P$ , $\lambda \in [0, \infty)$,
there exists $y\in P$ such that $x\perp y$ and $x+y\perp \lambda x-y$.
	
	An ordered pair $(X,\perp)$ is called an orthogonality space.
}

In 2010, Fechner and Sikorska \cite{23} studied the stability of orthogonality and proposed the definition of orthogonality as follows.

{\bf Definition 1.4}\label{dt2}
{\sl
	Let $X$ be an Abelian group and let $\perp$ be a binary relation defined on $X$ with the properties:

{\rm (1)} if $x,y\in X$ and $x\perp y$, then $x\perp -y,-x\perp y$ and $2x\perp 2y$;

{\rm (2)} for every $x\in X$, there  exists a $y \in X$ such that $x\perp y$ and $x+y\perp x-y$.
}

It's worth noting that every orthogonal space satisfies these conditions   as well as any normed linear space with the isosceles orthogonality, but Pythagorean orthogonality no longer satisfies these conditions.

Considering the current gaps, in this paper, we have made an attempt to prove the stability of orthogonal additivity in $\beta$-homogeneous  $F$-spaces and quasi-Banach spaces.

During the entire course of this work,  $\beta_2$ are considered as positive real numbers with  $\beta_2\geq0$. Furthermore, $X$ is assumed as an Abelian group while $Y$ is a $\beta_2$-homogeneous $F$-space.

\SUB{2. Stability of the orthogonally additive functional
equation}

\quad 
In this section, following some ideas from \cite{23}, we deal with the stability problem for the orthogonally additive functional equation in $\beta$-homogeneous  $F$-spaces. Our main theorem is the following.

{\bf Theorem 2.1} 
{\sl Let $X$ be an Abelian group, and $Y$ be a $\beta_2$-homogeneous $F$-space. 
	For $\varepsilon\geq 0$, assume $f: X\rightarrow Y$ be a mapping such that for all 
	$x, y \in X$ one has
\begin{align}\label{eqn1}
	x \perp y \quad \text { implies } \quad\|f(x+y)-f(x)-f(y)\| \leq \varepsilon. \tag {2.1}
	\end{align}
	Then there exists a mapping $g: X\rightarrow Y$ such that 
	\begin{align}\label{eqn2}
x \perp y \quad \text { implies } \quad g(x+y)=g(x)+g(y)
\tag {2.2}
	\end{align}
	and
	\begin{align}\label{eqn3}
\Vert f(x)-g(x)\Vert\leq b\varepsilon
	\tag {2.3}
	\end{align}
	with $$	b=\frac{2^{\beta_2+2}+3^{\beta_2+1}+3}{8^{\beta_2}}\left(1+\displaystyle\sum_{n=2}^{\infty}\left(\left(\frac{2^{n-1}+1}{2\cdot 4^{n-1}}\right)^{\beta_2}+\left(\frac{2^{n-1}-1}{2\cdot 4^{n-1}}\right)^{\beta_2}\right)\right)	$$
	for all $x \in 2X = \{2x : x \in X\}$. Moreover, the mapping $g$ is unique on the
	set $2X$.}

	{\bf Proof} By Definition 1.4, we get immediately there exists a $y \in X$ such 
that $x \perp y$ and $x+y \perp x-y$. Moreover, we also have 
$\pm 2x \perp \pm 2y, \pm (x+y) \perp \pm(x-y)$. Thus we conclude that
	$$
	\begin{aligned}\| 3 f(4 x)-8 f(2 x)-f(-4 x)) \| \leq 3^{\beta_2} &\|f(4 x)-f(2 x+2 y)-f(2 x-2 y)\| \\ &+\|f(-2 x+2 y)+f(-2 x-2 y)-f(-4 x)\| \\ &+3^{\beta_2}\|f(2 x+2 y)-f(2 x)-f(2 y)\| \\ &+3^{\beta_2}\|f(2 x-2 y)-f(2 x)-f(-2 y)\| \\ &+\|f(-2 x)+f(2 y)-f(-2 x+2 y)\| \\ &+\|f(-2 x)+f(-2 y)-f(-2 x-2 y)\| \\ &+2^{\beta_2}\|f(2 y)-f(y-x)-f(y+x)\| \\ &+2^{\beta_2}\|f(-2 y)-f(-y-x)-f(-y+x)\| \\ &+2^{\beta_2}\|f(y-x)+f(-y-x)-f(-2 x)\| \\ &+2^{\beta_2}\|f(y+x)+f(-y+x)-f(2 x)\| \\&\leq  (2^{\beta_2+2}+3^{\beta_2+1}+3)\varepsilon
	.\end{aligned}
	$$
	This proves that
	\begin{align}\label{eqn4}
\left\|f(2 x)-\frac{3}{8} f(4 x)+\frac{1}{8} f(-4 x)\right\| \leq \frac{2^{\beta_2+2}+3^{\beta_2+1}+3}{8^{\beta_2}}\varepsilon, \quad x \in X. \tag {2.4}
	\end{align}
	From now on, we set $a=\frac{2^{\beta_2+2}+3^{\beta_2+1}+3}{8^{\beta_2}}$ for convenience, then we have
	\begin{align}\label{eqn5}
\left\|f(2 x)-\frac{3}{8} f(4 x)+\frac{1}{8} f(-4 x)\right\| \leq a\varepsilon, \quad x \in X.\tag {2.5}
	\end{align}

	Now we will prove that for all $n \in \mathbb{N}$,
	\begin{align}\label{eqn6}
\left\|f(2 x)-\frac{2^{n}+1}{2 \cdot 4^{n}} f\left(2^{n+1} x\right)+\frac{2^{n}-1}{2 \cdot 4^{n}} f\left(-2^{n+1} x\right)\right\| \leq b \varepsilon, \quad x \in X, \tag {2.6}
	\end{align}
	
	with $	b=a\left(1+\displaystyle\sum_{n=2}^{\infty}\left(\left(\frac{2^{n-1}+1}{2\cdot 4^{n-1}}\right)^{\beta_2}+\left(\frac{2^{n-1}-1}{2\cdot 4^{n-1}}\right)^{\beta_2}\right)\right)	$.
	
	First, using (2.5), through a simple estimate we obtain
	$$
	\begin{aligned}&\left\|f(2 x)-\frac{2^{n+1}+1}{2 \cdot 4^{n+1}} f\left(2^{n+2} x\right)+\frac{2^{n+1}-1}{2 \cdot 4^{n+1}} f\left(-2^{n+2} x\right)\right\| \\& \leq\left\|f(2 x)-\frac{2^{n}+1}{2 \cdot 4^{n}} f\left(2^{n+1} x\right)+\frac{2^{n}-1}{2 \cdot 4^{n}} f\left(-2^{n+1} x\right)\right\| \\& +\bigg(\frac{2^{n}+1}{2 \cdot 4^{n}}\bigg)^{\beta_2}\left\|f\left(2^{n+1} x\right)-\frac{3}{8} f\left(2^{n+2} x\right)+\frac{1}{8} f\left(-2^{n+2} x\right)\right\| \\& +\bigg(\frac{2^{n}-1}{2 \cdot 4^{n}}\bigg)^{\beta_2}\left\|f\left(-2^{n+1} x\right)-\frac{3}{8} f\left(-2^{n+2} x\right)+\frac{1}{8} f\left(2^{n+2} x\right)\right\| \\& \leq\left\|f(2 x)-\frac{2^{n}+1}{2 \cdot 4^{n}} f\left(2^{n+1} x\right)+\frac{2^{n}-1}{2 \cdot 4^{n}} f\left(-2^{n+1} x\right)\right\|\\
	&+\left(\bigg(\frac{2^{n}+1}{2 \cdot 4^{n}}\bigg)^{\beta_2} +\bigg(\frac{2^{n}-1}{2 \cdot 4^{n}}\bigg)^{\beta_2} \right) a\varepsilon,\end{aligned}
    $$
    which implies
    $$
    \begin{aligned}&\left\|f(2 x)-\frac{2^{n+1}+1}{2 \cdot 4^{n+1}} f\left(2^{n+2} x\right)+\frac{2^{n+1}-1}{2 \cdot 4^{n+1}} f\left(-2^{n+2} x\right)\right\| \\& -\left\|f(2 x)-\frac{2^{n}+1}{2 \cdot 4^{n}} f\left(2^{n+1} x\right)+\frac{2^{n}-1}{2 \cdot 4^{n}} f\left(-2^{n+1} x\right)\right\|\\
    &\leq\left(\bigg(\frac{2^{n}+1}{2 \cdot 4^{n}}\bigg)^{\beta_2} +\bigg(\frac{2^{n}-1}{2 \cdot 4^{n}}\bigg)^{\beta_2} \right) a\varepsilon.\end{aligned}
    $$
    Now we let
    $$
    h(x,n)=\left\|f(2 x)-\frac{2^{n}+1}{2 \cdot 4^{n}} f\left(2^{n+1} x\right)+\frac{2^{n}-1}{2 \cdot 4^{n}} f\left(-2^{n+1} x\right)\right\|,
    $$
    so we have that
    $$
    h(x,n+1)-h(x,n)\leq\left(\bigg(\frac{2^{n}+1}{2 \cdot 4^{n}}\bigg)^{\beta_2} +\bigg(\frac{2^{n}-1}{2 \cdot 4^{n}}\bigg)^{\beta_2} \right) a\varepsilon,
    $$
    and then
    \begin{align*}
    h(x,n)&=\sum_{i=2}^{n}(h(x,i)-h(x,i-1))+h(x,1)\\
        &\leq\sum_{i=2}^{n}\left(\bigg(\frac{2^{i-1}+1}{2 \cdot 4^{i-1}}\bigg)^{\beta_2} +\bigg(\frac{2^{i-1}-1}{2 \cdot 4^{i-1}}\bigg)^{\beta_2} \right) a\varepsilon+a\varepsilon\\
        &\leq b\varepsilon,
    \end{align*}
    with
     $$
     b=\sum_{i=2}^{n}\left(\bigg(\frac{2^{i-1}+1}{2 \cdot 4^{i-1}}\bigg)^{\beta_2} +\bigg(\frac{2^{i-1}-1}{2 \cdot 4^{i-1}}\bigg)^{\beta_2}+1\right) a.
     $$
    This means that
    $$
    \left\|f(2 x)-\frac{2^{n}+1}{2 \cdot 4^{n}} f\left(2^{n+1} x\right)+\frac{2^{n}-1}{2 \cdot 4^{n}} f\left(-2^{n+1} x\right)\right\|\leq b\varepsilon,\quad x\in X.
    $$
    	The next step is to prove that for each $x \in X$ the sequence
	$$
	g_{n}(x):=\frac{2^{n}+1}{2 \cdot 4^{n}} f\left(2^{n} x\right)-\frac{2^{n}-1}{2 \cdot 4^{n}} f\left(-2^{n} x\right), \quad n \in \mathbb{N}
	$$
	is convergent in $Y$. Since $Y$ is complete, it suffices to show that $(g_n(x))_{n\in \mathbb{N}}$ is 
	a Cauchy sequence for every $x \in X$. Applying estimate (2.5) twice then we have
	$$
	\begin{aligned}\left\|g_{n}(x)-g_{n+1}(x)\right\|=& \bigg\| \frac{2^{n}+1}{2 \cdot 4^{n}}\left(f\left(2^{n} x\right)-\frac{3}{8} f\left(2^{n+1} x\right)+\frac{1}{8} f\left(-2^{n+1} x\right)\right) \\ &-\frac{2^{n}-1}{2 \cdot 4^{n}}\left(f\left(-2^{n} x\right)-\frac{3}{8} f\left(-2^{n+1} x\right)+\frac{1}{8} f\left(2^{n+1} x\right)\right) \bigg\| \\ 
	\leq & \left(\bigg(\frac{2^{n}+1}{2 \cdot 4^{n}}\bigg)^{\beta_2}+\bigg(\frac{2^{n}-1}{2 \cdot 4^{n}}\bigg)^{\beta_2}\right) a \varepsilon \end{aligned}
	$$
	for each $n \in \mathbb{N}$, which gives us  that $(g_n(x))_{n\in N}$ is a Cauchy sequence.
	
 Hence, the mapping $g : X \rightarrow Y$ can be defined as:
    $$ 
	g(x):=\lim _{n \rightarrow \infty} g_{n}(x)
	$$
	for all $x \in X$.
	Combining with (2.6) we have
	$$\Vert f(2x) - g(2x)\Vert \leq  b\varepsilon, \quad x \in X.$$
	In order to prove that $g$ is orthogonally additive observe first that for
	$x, y \in X$ such that $x \perp y$ and $n \in N, n > 1$ we have
	$$
	\begin{aligned}&\quad\left\|g_{n}(x+y)-g_{n}(x)-g_{n}(y)\right\| \\
	&= \bigg\| \frac{2^{n}+1}{2 \cdot 4^{n}} f\left(2^{n}(x+y)\right)-\frac{2^{n}-1}{2 \cdot 4^{n}} f\left(-2^{n}(x+y)\right) \\ 
	&\quad-\frac{2^{n}+1}{2 \cdot 4^{n}} f\left(2^{n} x\right)+\frac{2^{n}-1}{2 \cdot 4^{n}} f\left(-2^{n} x\right)-\frac{2^{n}+1}{2 \cdot 4^{n}} f\left(2^{n} y\right)+\frac{2^{n}-1}{2 \cdot 4^{n}} f\left(-2^{n} y\right) \bigg\| \\
	&=\bigg \| \frac{2^{n}+1}{2 \cdot 4^{n}}\left[f\left(2^{n}(x+y)\right)-f\left(2^{n} x\right)-f\left(2^{n} y\right)\right] \\ 
	&\quad-\frac{2^{n}-1}{2 \cdot 4^{n}}\left[f\left(2^{n}(-x-y)\right)-f\left(-2^{n} x\right)-f\left(-2^{n} y\right)\right]\bigg \| \\
	& \leq  \left(\frac{2^{n}+1}{2 \cdot 4^{n}}\right)^{\beta_2}\left\|f\left(2^{n}(x+y)\right)-f\left(2^{n} x\right)-f\left(2^{n} y\right)\right\| \\ 
	&\quad+\left(\frac{2^{n}-1}{2 \cdot 4^{n}}\right)^{\beta_2}\left\|f\left(2^{n}(-x-y)\right)-f\left(-2^{n} x\right)-f\left(-2^{n} y\right)\right\| \\
	& \leq  \left(\left(\frac{2^{n}-1}{2 \cdot 4^{n}}\right)^{\beta_2}+\left(\frac{2^{n}-1}{2 \cdot 4^{n}}\right)^{\beta_2}\right) \varepsilon.
	\end{aligned}
	$$
	 Moreover, letting $n \rightarrow \infty$,  we get (2.2).

	Now, we show the uniqueness of $g$. Assuming $g^{\prime}$ as another   mapping satisfying 
	(2.2) and (2.3)  that yields:

	$$
	\left\|g(x)-g^{\prime}(x)\right\| \leq\|g(x)-f(x)\|+\left\|g^{\prime}(x)-f(x)\right\| \leq 2b \varepsilon
	$$
	for all $x\in 2X$.

	On the other hand, the mapping $g-g^{\prime}$ satisfies (2.2) and thus, in
	particular, (2.1) with $\varepsilon= 0$. By applying (2.6) to $g-g^{\prime}$ we see that
	$$
	\begin{aligned} g(2 x)-g^{\prime}(2 x)=& \frac{2^{n}+1}{2 \cdot 4^{n}}\left[g\left(2^{n+1} x\right)-g^{\prime}\left(2^{n+1} x\right)\right] \\ &-\frac{2^{n}-1}{2 \cdot 4^{n}}\left[g\left(-2^{n+1} x\right)-g^{\prime}\left(-2^{n+1} x\right)\right] \end{aligned}
	$$
	and  therefore
	$$
	\begin{aligned}\left\|g(2 x)-g^{\prime}(2 x)\right\| \leq & \left(\frac{2^{n}+1}{2 \cdot 4^{n}}\right)^{\beta_2}\left\|g\left(2^{n+1} x\right)-g^{\prime}\left(2^{n+1} x\right)\right\| \\ &+\left(\frac{2^{n}-1}{2 \cdot 4^{n}}\right)^{\beta_2}\left\|g\left(-2^{n+1} x\right)-g^{\prime}\left(-2^{n+1} x\right)\right\| \\ 
	\leq & \left(\left(\frac{2^{n}-1}{2 \cdot 4^{n}}\right)^{\beta_2}+\left(\frac{2^{n}-1}{2 \cdot 4^{n}}\right)^{\beta_2}\right) 2b \varepsilon
	\end{aligned}
	$$
	for $x \in X$. 
	
	Combining the both inequalities, we  can easily get the thesis.
$\Box$

	{\bf Remark 2.2}
	{\sl
		Compare with the thereom of \cite{23}, $Y$ is a $F$-space with $\beta_2$-homogenous instead of homogenous of \cite{23}, so the estimate combined with $\beta_2$-homogenuous of (2.6) here would be more difficult and complicated and moreover we can get (2.6) is small enough so that we can get the convergence of the function sequence $(g_n(x))_{n\in N}$,that is the main point here.
	}

By the same method, we can also obtain the stability result for different target spaces, where the  space $Y$  is equipped with quasi-norm. From now on, in corollaries, assume that $X$ is an Abelian group and $Y$ 
is a quasi-Banach space.

{\bf Corollary 2.3}
{\sl
For $\varepsilon\geq 0$, let $f: X\rightarrow Y$ be a mapping such that for all 
$x, y \in X$ one has
$$
x \perp y \quad \text { implies } \quad\|f(x+y)-f(x)-f(y)\| \leq \varepsilon
$$
Then there exists a mapping $g: X\rightarrow Y$ such that 
$$
x \perp y \quad \text { implies } \quad g(x+y)=g(x)+g(y)
$$
and
$$\Vert f(x)-g(x)\Vert\leq b^{\frac{1}{p}}\varepsilon$$
with $$	b=\frac{2^{p+2}+3^{p+1}+3}{8^{p}}\left(1+\displaystyle\sum_{n=2}^{\infty}\left(\left(\frac{2^{n-1}+1}{2\cdot 4^{n-1}}\right)^{p}+\left(\frac{2^{n-1}-1}{2\cdot 4^{n-1}}\right)^{p}\right)\right)	$$
for all $x \in 2X = \{2x : x \in X\}$. Moreover, the mapping $g$ is unique on the
set $2X$.}

{\bf Proof}
	Let $\|\cdot\|_p=\|\cdot\|^p$, then it is obviously that $(Y, \|\cdot\|_p)$ is $p$-homogeneous, we obtain
$$
x \perp y \quad \text { implies } \quad\|f(x+y)-f(x)-f(y)\|_p \leq \varepsilon^p.  	
$$
According to Theorem 2.1, we obtain that there exists a mapping $g: X\rightarrow Y$ such that 
$$
x \perp y \quad \text { implies } \quad g(x+y)=g(x)+g(y)
$$
and
$$\Vert f(x)-g(x)\Vert_p\leq b\varepsilon^p$$
with $	b=\frac{2^{p+2}+3^{p+1}+3}{8^{p}}\left(1+\displaystyle\sum_{n=2}^{\infty}\left(\left(\frac{2^{n-1}+1}{2\cdot 4^{n-1}}\right)^{p}+\left(\frac{2^{n-1}-1}{2\cdot 4^{n-1}}\right)^{p}\right)\right)	$
for all $x \in 2X = \{2x : x \in X\}$. Moreover, the mapping $g$ is unique on the
set $2X$ and the claim follows.
$\Box$
\mbox{}
\\[8pt]
{\bf Acknowledgements}\ \ We thank the referees for their time and comments.

\end{document}